\documentclass[12pt]{article}
\usepackage{amsmath,amsfonts,amsthm}
\begin {document}

\title{Simple subalgebras of simple special Jordan algebras}
\author{
M.V.Tvalavadze\thanks{\it Email address: marina@math.mun.ca}, \thanks{supported by NSERC Grant 227060-00}\\
 \small Department of Mathematics and Statistics\\
\small Memorial University of Newfoundland\\
\small St. John's, NL, CANADA. }
\date{}
\maketitle

\newtheorem{theorem}{Theorem}[section]

\newtheorem{corollary}[theorem]{Corollary}
\newtheorem{lemma}[theorem]{Lemma}

\newtheorem{definition}[theorem]{Definition}

\sloppy

\begin{abstract}
In this paper we determine  all types and the canonical forms of
simple subalgebras for each type of simple Jordan algebras and the
number of conjugate classes corresponding to the given simple
Jordan algebra.
\end{abstract}

\section{Introduction}
This paper provides a classification of simple subalgebras in
finite-dimensional special simple Jordan algebras over
algebraically closed field $F$ with characteristic unequal to  2.
More precisely, we will determine the canonical forms for any
simple subalgebras of special simple Jordan algebras, and the
number of conjugate classes corresponding to the given simple
Jordan subalgebra. In particular, a Jordan algebra of any type can
be realized as a Jordan subalgebra of symmetric or symplectic
matrices of an appropriate order.

In 1987 N. Jacobson determined the orbits under the orthogonal
group $O(n)$ of the subalgebras of the Jordan algebra of $n\times
n$ real symmetric matrices (\cite{jac}).

The paper significantly relies on the description of maximal
subalgebras of finite-dimensional special simple Jordan algebras
obtained by M. Racine in 1974 (\cite{R2}).

In this paper we consider simple Jordan algebras presented in
canonical matrix realizations, that is, $H(F_n)$ can be viewed as
the algebra of all symmetric matrices in $F_n$ with respect to the
ordinary transpose;  $F^{(+)}_n$ is the set of all matrices of
order $n$ under the circle product $A\odot B=\frac{AB+BA}{2}$;
$H(F_{2n},j)$ where $j$ denotes a symplectic involution consists
of all matrices of order $2n$ of the form
$$        \left( \begin{array}{cc}
                   A&  B\\
                   C& A^t
                    \end{array}\right),
$$

\noindent where $B$, $C$ are any skew-symmetric matrices of order
$n$, and $A$ is any matrix of order $n$. If $f$ is a non-singular
symmetric bilinear form on a vector space $V$, then ${\cal
J}=F\oplus V$ is a Jordan algebra of the type $J(f,1)$.

Throughout the paper we assume that the base field $F$ is
algebraically closed with characteristic not two.

\section{Subalgebras}

\subsection{Matrix subalgebras}

Let $\cal J$ be a simple Jordan algebra of the type
$F_{\frac{n}{2}}^{(+)}$ where $n$ is even. Then it can always be
presented as a subalgebra of $H(F_n)$ as follows
$$  \left\{ \left(\begin{array}{cc}
            A& B\\
            -B& A
            \end{array}\right)\right\},\eqno(1)
$$
where $A$ is any symmetric matrix of order $\frac{n}{2}$ and $B$
is any skewsymmetric matrix of order $\frac{n}{2}$.

\begin{lemma}
Any automorphism of a Jordan algebra of the form (1) is induced by
an automorphism of $H(F_n)$.
\end{lemma}
\begin{proof}
Any automorphism of $\cal J$ can be extended to an automorphism or
antiautomorphism of a special universal enveloping algebra
$U({\cal J})$ which is isomorphic to $F_{\frac{n}{2}}\oplus
F_{\frac{n}{2}}$ (see \cite{jac2}). Notice that in this particular
case the associative enveloping algebra $S({\cal J})$ is
isomorphic to $U({\cal J})$ because from the explicit form (1)
$S({\cal J})$ consists of all matrices of the form:
$$  \left\{ \left(\begin{array}{cc}
            X& Y\\
           -Y&X
            \end{array}\right)\right\},
$$
where $X$ and $Y$ are any matrices of order $\frac{n}{2}$. Since
any automorphism of $F_{\frac{n}{2}}\oplus F_{\frac{n}{2}}$ either
induces non-trivial automorphisms of these ideals or sends one
ideal onto another, it can be lifted up to an inner automorphism
of the entire matrix algebra $F_n$. Consequently, for any
antiautomorphism of $F_{\frac{n}{2}}\oplus F_{\frac{n}{2}}$ we can
choose an automorphism (not necessarily non-trivial) of
$F_{\frac{n}{2}}\oplus F_{\frac{n}{2}}$ such that their
composition induces non-trivial antiautomorphisms of simple
ideals. Therefore, any (Jordan) automorphism of $\cal J$  can be
written as follows:
 $$ \varphi(X)=Q^{-1}XQ$$
or
$$\varphi(X)=Q^{-1}X^tQ, \eqno (2)$$

\noindent for some non-singular matrix $Q$.

 The next step is to prove that $\varphi$ is orthogonal.
In other words, all we have to show is that for any automorphism
$\varphi$ of $\cal J$, we can choose $Q$ such that (2) holds and
$Q^tQ=I$ where $I$ is the identity matrix. Since $\cal J$ is a
subalgebra of $H(F_n)$, for each $X$ in $\cal J$,
$(Q^{-1}XQ)^t=Q^{-1}XQ$, $Q^tX(Q^{-1})^t=Q^{-1}XQ$, $QQ^tX=XQQ^t$.
Denote $B=QQ^t$. Next we are going to show that $B$ is actually a
scalar multiple of the identity matrix. We are given that $BX=XB$
where $X$ is any matrix of the form (1). Let us write $B$ as
follows:
$$ B=\left(\begin{array}{cc}
            B_1& B_2\\
            B_3& B_4
            \end{array}\right)
$$
where $B_i$ are matrices of order $\frac{n}{2}$. By performing the
matrix multiplication, we obtain the $B_2=B_3=0$, and
$B_1=B_4=\alpha I$, for some non-zero $\alpha$. Since the ground
field $F$ is algebraically closed, we can choose $\beta\in F$ such
that $\alpha=\beta^2$. Set $Q'=\beta^{-1}Q$. Obviously, $Q'$
determines the same automorphism as $Q$ does, and $Q'^tQ'=I$.

The Lemma is proved.
\end{proof}

\begin{lemma} Let ${\cal J}\cong F_{\frac{n}{2}}^{(+)}$ be a
subalgebra of $H(F_n)$. Then, by an appropriate automorphism of
$H(F_n)$, $\cal J$ can always be reduced to the form (1).
\end{lemma}
\begin{proof}
Since $\cal J$ has the type $F_{\frac{n}{2}}^{(+)}$, by some (not
necessarily orthogonal) automorphism $\varphi$ of $F_n^{(+)}$  we
can always bring $\cal J$ to the following form (see \cite{jac2}
and \cite{tm})
$$ \left\{\left(\begin{array}{cc}
            X& 0\\
            0& X^t
            \end{array}\right)\right\}\eqno(3)
$$
where $X$ is any matrix of order $\frac{n}{2}$. Then, it is easily
seen that $\theta(Y)=S^{-1}YS$, where $S=\left(\begin{array}{cc}
            I            & iI\\
            \frac{1}{2}I  & -\frac{i}{2}I
            \end{array}\right)$, $I$ is the identity matrix, $i^2=-1$, sends
each element of the form (3) into the algebra of the form (1).
Therefore, by $\chi=\theta\circ\varphi$ we can bring $\cal J$ to
the form (1).

Next we will show that $\chi$ is actually an orthogonal
automorphism. Notice that $\chi$ sends $H(F_n)$ onto a Jordan
subalgebra of $F_n^{(+)}$ which consists of all matrices symmetric
with respect to the following involution: $j'=\chi\circ t\circ
\chi^{-1}$ where $t$ is the standard transpose involution. This
involution can be rewritten as follows $j'(X)=C^{-1}X^tC$ for some
non-singular symmetric matrix $C$ of order $n$. It follows from
the above considerations that any matrix of the form (1) is
symmetric with respect to $j'$. Equivalently, for any $Y$ of the
form (1), $C^{-1}Y^tC=Y$, $Y^tC=CY$, $YC=CY$ because $Y$ is
symmetric. As proved in the previous Lemma, $C=\alpha I$ for some
non-zero $\alpha$. Therefore, $j'=t$, and $\chi(H(F_n))=H(F_n)$,
and $\chi$ is actually an automorphism of $H(F_n)$. Hence, the
Lemma is proved.
\end{proof}

\begin{lemma}
Let $\cal A$ be a special simple matrix Jordan algebra, and $\cal
J$ be a proper simple subalgebra of $\cal A$. Denote a maximal
subalgebra which contains $\cal J$ as $M$. Next, consider a
Wedderburn splitting $M=S\oplus R$ where $S$ is a semisimple
algebra, $R$ is the radical. Then, there exists an automorphism
$\varphi$ of $\cal A$ such that $\varphi({\cal J})\subseteq S$.
\end{lemma}
\begin{proof}
Let 1 be the identity element of $\cal A$, and $1\in {\cal J}$.
According to \cite{Mc2}, if $\cal J$ is special and the degree of
$\cal J$ is not divisible by the characteristic, then $\cal J$ is
conjugate under an inner automorphism $T$ of $M$ to some
subalgebra of $S$, and  $T$ is a composition of the standard
automorphisms $T_{x,y}$ that can be represented  in associative
terms as follows
$$
T_{x,y}(a)=tat^{-1}, \eqno(4)$$ where
 $t=u^{-\frac{1}{2}}(1-xy)(1+yx)$, $u=(1-xy)(1+yx)(1+xy)(1-yx),$
$x,y\in M$.

Let $x$, $y$ be symmetric with respect to an involution $j$ of
$\cal A$: $j(x)=x$, $j(y)=y$. Then, it is obvious that

 $$ j(u)=u, \quad \mbox{and} \quad j(t)=t^{-1}. \eqno(5)$$

 If ${\cal A}=F_n^{(+)}$, then from the explicit form (4)
$T_{x,y}$ is easily extendable to $\cal A$. If ${\cal A}=H(F_n)$,
then, because of (5), $T_{x,y}$ is orthogonal, therefore,
extendable to $\cal A$. If ${\cal A}=H(F_{2n},j)$, then, because
of (5),  $T_{x,y}$ is symplectic, therefore, extendable to $\cal
A$.

If  $\cal J$ is special and the degree of $\cal J$ is divisible by
characteristic, then $T$ is a {\it generalized} inner automorphism
(see \cite{Mc2}), that is, $T$ is a composition of an
automorphisms $T_{x_1,\ldots,x_n,m}$ of the form
$$T_{x_1,\ldots,x_n,m}=U_v^{-1}(I+V_{x_1,\ldots,x_n,m}+U_{x_1}\ldots U_{x_n}U_{-m})
(I+V_{m,x_n\ldots,x_1}+U_{m}U_{x_n}\ldots U_{x_1})$$

\noindent where $v, x_i\in M$, $m\in R$. In associative terms
operators take the form:
$$U_v(a)=vav,$$
$$U_{x_i}(a)=x_iax_i,$$
$$V_{x_1,\ldots,x_n,m}(a)=x_1\ldots x_nma+amx_n\ldots x_1.$$

Hence, if all $x_i$, $m$ and $a$ are symmetric with respect to an
involution of $\cal A$, then
$j(T_{x_1,\ldots,x_n,m}(a))=T_{x_1,\ldots,x_n,m}(a)$. Therefore,
$T_{x_1,\ldots,x_n,m}$ as well as $T$ is extendable to $\cal A$.
The Lemma is proved.
\end{proof}

\begin{lemma} Let ${\cal J}\cong F_n^{(+)}$ be a
subalgebra of $H(F_{2n},j)$. Then, by an appropriate automorphism
of $H(F_{2n},j)$, $\cal J$ can always be reduced to the following
form
$$\left\{ \left(\begin{array}{cc}
            X& 0\\
            0& X^t
            \end{array}\right)\right\}\eqno(6)
$$

\end{lemma}
\begin{proof}
Since $\cal J$ has the type $F_n^{(+)}$, by some  automorphism
$\varphi$ of $F_{2n}^{(+)}$, $\cal J$ can be brought  to the form
(6) (see \cite{jac2} and \cite{tm}). Notice that $\varphi$ sends
$H(F_{2n},j)$ onto a Jordan subalgebra of $F_{2n}^{(+)}$ which
consists of all matrices symmetric with respect to the following
involution: $j'=\varphi\circ j\circ \varphi^{-1}$. This involution
can be rewritten as follows $j'(Y)=C^{-1}Y^tC$ for some
non-singular skew-symmetric matrix $C$ of order $2n$. It follows
from the above considerations that any matrix of the form (6) is
symmetric with respect to $j'$. Equivalently, for any $Y$ of the
form (6), $C^{-1}Y^tC=Y$, $Y^tC=CY$. Acting in the same manner as
above, we can show that $C=\alpha \left(\begin{array}{cc}
            0& I_n\\
            -I_n& 0
            \end{array}\right) $
for some non-zero $\alpha$, where $I_n$ denotes the identity
matrix of order $n$. Therefore, $\varphi(H(F_n,j))=H(F_n,j)$, and
$\varphi$ is an automorphism of $H(F_n,j)$. Hence, the Lemma is
proved.
\end{proof}

\begin{definition}
Subalgebras ${\cal J}_1$ and ${\cal J}_2$ of a Jordan algebra
$\cal A$ are said to be equivalent if there exists an automorphism
$\varphi$ of $\cal A$ such that  either ${\cal J}_1=\varphi({\cal
J}_2)$ or ${\cal J}_2=\varphi({\cal J}_1)$.
\end{definition}

\begin{definition}
Let $\cal J$ be a subalgebra of $\cal A$. Then the set $C({\cal
J})$ of all subalgebras equivalent to $\cal J$ in $\cal A$ is said
to be a conjugate class of $\cal J$.
\end{definition}
\par\medskip
\par\medskip
{\bf\large Canonical realizations of simple subalgebras}
\par\medskip
Let $\cal A$ be a simple Jordan algebra, and $\cal J$ be a simple
subalgebra of $\cal A$. All realizations listed below we will call
{\it canonical}.
\par\medskip
\noindent {\bf 1}. Let ${\cal A}=F_n^{(+)}$
\par\medskip
\noindent (1.1)  ${\cal J}\cong F_m^{(+)}$, ${\cal
J}=\{\mbox{diag}(X,\ldots,X,X^t,\ldots,X^t,0,\ldots,0)\}$ where
$X$ is any matrix of order $m$.

\noindent (1.2) ${\cal J}\cong H(F_m)$, ${\cal
J}=\{\mbox{diag}(X,\ldots,X, 0,\ldots,0)\}$ where $X$ is any
symmetric matrix of order $m$.

\noindent (1.3) ${\cal J}\cong H(F_{2m},j)$, ${\cal
J}=\{\mbox{diag}(X,\ldots,X,0,\ldots,0)\}$ where $X$ is any
symplectic matrix of order $2m$.

\par\medskip
\noindent {\bf 2}. Let  ${\cal A}=H(F_n)$
\par\medskip
\noindent (2.1) ${\cal J}\cong F_m^{(+)}$, ${\cal
J}=\{\mbox{diag}(X,\ldots,X,0\ldots,0)\}$ where $X$ is of the form
(1) in which $A$ and $B$ are of order $m$.

\noindent (2.2) ${\cal J}\cong H(F_m)$, ${\cal
J}=\{\mbox{diag}(X,\ldots,X,0,\ldots,0)\}$ where $X$ is any
symmetric matrix of order $m$.

\noindent (2.3) ${\cal J}\cong H(F_{2m},j)$, ${\cal
J}=\{\mbox{diag}(X,\ldots,X,0\ldots,0)\}$
$$X=\left(\begin{array}{cccc}
             A& -B& -C& D \\
             B&  A&  D& C \\
             C&  D&  A&-B \\
             D& -C&  B& A
\end{array}\right)$$
where $A$ is a symmetric matrix of order $m$, $B$, $C$, $D$ are
skew-symmetric matrices of order $m$.

\par\medskip
\noindent {\bf 3}. Let ${\cal A}=H(F_{2n},j)$
\par\medskip
\noindent (3.1)  ${\cal J}\cong F_m^{(+)}$, $${\cal
J}=\{\mbox{diag}(\underbrace{X,\ldots,X,}_{k}\underbrace{X^t,\ldots,X^t,}_{l}
\underbrace{0,\ldots,0}_{s},X^t,\ldots,X^t,X,\ldots,X,0,\ldots,0)\}$$
where $k+l+s=n$,  $X$ is any matrix of order $m$.

\noindent (3.2) ${\cal J}\cong H(F_m)$, ${\cal
J}=\{\mbox{diag}(\underbrace{X,\ldots,X,}_{k}\underbrace{0,\ldots,0}_{l},
X,\ldots,X,0,\ldots,0)\}$ where $k+l=n$, $X$ is any symmetric
matrix of order $m$.

\noindent (3.3) ${\cal J}\cong H(F_{2m},j)$, ${\cal J}=\left\{
\left(
\begin{array}{cc}
                   A&  B\\
                   C& A^t
                    \end{array}\right)\right\}$
                    $$ A=\small \left(\begin{array}{cccccccc}
  X&        &   &   &    &      &    &    \\
   & \ddots &   &   &    &   0  &    &   \\
    &       & X &   &    &      &    &    \\
    &       &   &  X&   Y&      &    &    \\
    &       &   &  Z& X^t&      &    &    \\
    &       & 0 &   &    &\ddots&    &   \\
    &       &   &   &    &      &   X& Y   \\
    &       &   &   &    &      &   Z& X^t
         \end{array}\right),
$$

$$ B=\small \left(\begin{array}{cccccccc}
  Y&        &   &   &    &      &    &    \\
   & \ddots &   &   &    &    0  &    &   \\
    &       & Y &   &    &      &    &    \\
    &       &   &   &    &      &    &    \\
    &       &   &   &    &      &    &    \\
    &       &  0 &   &    &  0  &    &   \\
    &       &   &   &    &      &    &    \\
    &       &   &   &    &      &    &
         \end{array}\right),
   C=\small \left(\begin{array}{cccccccc}
  Z&        &   &   &    &      &    &    \\
   & \ddots &   &   &    &    0  &    &   \\
    &       & Z &   &    &      &    &    \\
    &       &   &   &    &      &    &    \\
    &       &   &   &    &      &    &    \\
    &       &  0 &   &    &  0  &    &   \\
    &       &   &   &    &      &    &    \\
    &       &   &   &    &      &    &
         \end{array}\right)
$$
where $X$ is any matrix of order $m$, $Y$,$Z$ are skew-symmetric
matrices of order $m$.
\begin{definition}
Let $\cal J$ and ${\cal J}'$ be two proper subalgebras of $\cal
A$, and ${\cal J}'$ be given in the canonical realization. If
$\cal J$ is equivalent to ${\cal J}'$, then ${\cal J}'$ is said to
be the {\it canonical form} of ${\cal J}$.
\end{definition}

\begin{theorem}Let $\cal A$ be a simple matrix Jordan algebra. Then,
any simple matrix subalgebra of $\cal A$ has a unique canonical
form as above.
\end{theorem}
\begin{proof}
Let $\cal J$ be any proper simple matrix subalgebra of $\cal A$.
In particular,  the degree of ${\cal J}\ge 3$. Denote the identity
of $\cal A$ as 1.

The proof of the Theorem consists of three cases.

{\bf Case 1}  ${\cal A}= F_n^{(+)}$

{\bf 1.1} Let $\cal J$ be of the type $F_m^{(+)}$ for some $m<n$.
Since any Jordan algebra of this type has precisely two
non-equivalent irreducible representations both of which have
degree $m$ (see \cite{jac2}), $\cal J$ is equivalent to the
subalgebra in the canonical realization (1.1). If $1\in {\cal J}$,
then the last zero in (1.1) is omitted. Next we are going to show
that $\cal J$ has a unique canonical form. Equivalently, if $l$
and $k$ are the number of $X$-blocks and $X^t$-blocks in (1.1),
then $|l-k|$ is an invariant for $\cal J$. Indeed, let $S({\cal
J})$ be a simple algebra. Then, either $l$ or $k$ is zero, and
$|l-k|=\mbox{rk}\,e$ where $e$ is the identity of $\cal J$. If
$S({\cal J})$ is a non-simple semisimple algebra, then $S({\cal
J})={\cal I}_1\oplus {\cal I}_2$ where ${\cal I}_i$ are simple
ideals with the identity elements $e_i$. Hence,
$|l-k|=|\mbox{rk}\,e_1-\mbox{rk}\,e_2|$ which is invariant for
$\cal J$. Denote $k_{\cal A}({\cal J})=|l-k|$.

{\bf 1.2} Let $\cal J$ be of the type $H(F_m)$ for some $m\le n$.
Then, it follows from the uniqueness of the irreducible
representation of  $H(F_m)$ (see \cite{jac2}) that  $\cal J$ is
equivalent to the subalgebra in the canonical realization (1.2).
If $1\in {\cal J}$, then the last zero in (1.2) is omitted.

{\bf 1.3}  The proof of the case when ${\cal J}\cong H(F_{2m},j)$,
$2m<n$, is exactly the same as the previous proof. In particular,
$\cal J$ of the type $H(F_{2m},j)$ is equivalent to the subalgebra
in the canonical realization (1.3). Obviously, the canonical form
is unique.
\par\medskip
{\bf Case 2}  ${\cal A}=H(F_n)$

Here, our main goal is to determine the canonical form of any
simple matrix Jordan subalgebra of $H(F_n)$.  Let $M$ be a maximal
subalgebra of $H(F_n)$. According to \cite{R2}, $M$  is isomorphic
to one of the following:

1. $H(F_k)\oplus H(F_l)$, $k+l=n$,

2. $F_k^{(+)}\oplus H(F_l)\oplus R$, $2k+l=n$, $R$ is the radical
(if $l=0$, then $M\cong F_{\frac{n}{2}}^{(+)}\oplus R$)

3. $J(f,1)$ only if $n=2^{m}$ and either $\dim\,J(f,1)=2(m+1)$,
$m$ is even, or $\dim\,J(f,1)=2m+1$, $m$ is odd.

First, assume that $\cal J$ is a simple matrix subalgebra of
$H(F_n)$ such that $1\in {\cal J}$. Then there exists a maximal
subalgebra $M$ such that ${\cal J}\subset M$. Since
$\mbox{deg}\,{\cal J}\ge 3$, $M$ cannot be of the type 3.
 If $M$ contains a non-zero radical, that is, $M=S\oplus R$, where $S$ a semisimple algebra,
 $R$ the radical, then by Lemma 2.3 we
can  assume that ${\cal J}\subseteq S$. If $S=S_1\oplus S_2$ where
$S_i$ non-trivial simple ideals, we can choose three orthogonal
idempotents (see \cite{R2}): $e$, $e^t$, $ff^t$, $1=e+e^t+ff^t$
such that
$$ S_1=ff^tH(F_n)ff^t,\quad  S_2=eF_ne+e^tF_n^te^t\eqno (7)$$

 Since $ff^t$ is an element of $H(F_n)$, by an automorphism $\varphi$ of
$H(F_n)$, it can be reduced to the  following form:
$$ \varphi(ff^t)= \left(\begin{array}{cc}
            I_l& 0\\
            0& 0
            \end{array}\right)
$$
where $I_l$ is the identity matrix of order $l$. Since $e$ and
$e^t$ are orthogonal to $ff^t$, they take the forms:
$$\varphi(e)=\left(\begin{array}{cc}
            0& 0\\
            0& K
            \end{array}\right),
 \varphi(e^t)=\left(\begin{array}{cc}
            0& 0\\
            0& K^t
            \end{array}\right),
$$
where $K$ is a matrix of order $n-l$. Therefore, according to (7),
$$\varphi(S)= \left\{ \left(\begin{array}{cc}
            X& 0\\
            0& Y
            \end{array}\right)\right\},
  \varphi(S_1)= \left\{ \left(\begin{array}{cc}
            X& 0\\
            0& 0
            \end{array}\right)\right\},$$
$$\varphi(S_2)= \left\{ \left(\begin{array}{cc}
            0& 0\\
            0& Y
            \end{array}\right)\right\},\eqno(8)
$$
where $X$ is any symmetric matrix of order $l$, $Y$ is a symmetric
matrix of order $2k=n-l$ which is also an element of a subalgebra
of the type $F^{(+)}_k$.

In the case when $M$ is semisimple, that is, $M=S=S_1\oplus S_2$,
there exist two orthogonal idempotents such that
$$ S=eH(F_n)e+fH(F_n)f,\quad e+f=1.$$
Acting in the same manner as above we can reduce $S$ to (8).

Therefore, we can  define two homomorphisms $\pi_1$, $\pi_2$ as
projections on $S_1$ and $S_2$, respectively. Since $1\in {\cal
J}$, $\pi_1({\cal J})\ne\{0\}$,  $\pi_2({\cal J})\ne\{0\}$. This
implies that ${\cal J}\cong \pi_1({\cal J})\subseteq H(F_l)$,
$l<n$, and ${\cal J}\cong \pi_2({\cal J})\subseteq H(F_{2k})$,
$2k<n$. Therefore, we can reduce the problem of finding the
canonical form of $\cal J$ to the case of all  symmetric matrices
of order less than $n$. However, the above reduction does not work
in the case when $S$ is simple, that is,
$M=F_{\frac{r}{2}}^{(+)}\oplus R$, $r\le n$.
  Hence we can conclude that as soon as the given
simple subalgebra $\cal J$ is in the maximal subalgebra $M$ which
has a non-simple semisimple factor $S$, the problem can be reduced
to the case of symmetric matrices of a lower order. This process
stops only if at some step either $\pi_i({\cal J})\subseteq M\cong
F^{(+)}_{\frac{r}{2}}\oplus R$, or $\pi_i({\cal J})$ coincides
with $S_i$. Without any loss of generality,  we can assume that
$r=n$, that is,  ${\cal J}\subseteq M\cong
F^{(+)}_{\frac{n}{2}}\oplus R$.

All we need to reach our goal is to determine the canonical form
of $\cal J$  which is covered by a maximal subalgebra of the type
$F_{\frac{n}{2}}^{(+)}\oplus R$.   Notice that there is an
isomorphic imbedding $\theta$ of $F_{\frac{n}{2}}^{(+)}$ into
$H(F_n)$ such that $\theta(A+iB)=\left(\begin{array}{cc}
            A& B\\
            -B& A
            \end{array}\right)$, where $A$ is a symmetric matrix of order $\frac{n}{2}$,
$B$ is a skew-symmetric matrix of order $\frac{n}{2}$, $i^2=-1$.

{\bf 2.1}. Let us assume that $\cal J$ has the type $F^{(+)}_m$
where $n=2ml$. We know that by an appropriate automorphism $\psi$
of $F^{(+)}_{\frac{n}{2}}$, we can reduce $\theta^{-1}({\cal J})$
to the following canonical form:

$$\psi(\theta^{-1}({\cal J}))= \left\{ \left(\begin{array}{cccccc}
         X&  \ldots   &      0 &         0 &  \ldots &  0 \\
\vdots    &    \ddots &  \vdots&  \vdots   &  \ddots &  \vdots \\
   0      &  \ldots   &      X &         0 &  \ldots &  0  \\
   0      & \ldots    &      0 &         X^t &  \ldots &  0   \\
\vdots    & \ddots    & \vdots &   \vdots  &  \ddots &  \vdots  \\
   0      & \ldots    &      0 &         0 &  \ldots &       X^t
    \end{array}\right)\right\},
$$
where $X$ is any matrix of order $m$. Then, $X$ can be written as
$A+iB$ for an appropriate symmetric $A$ and skew-symmetric $B$.
Therefore, $\theta(\psi(\theta^{-1}({\cal J})))$ has  the
following representation  in $H(F_n)$:
$$\theta(\psi(\theta^{-1}({\cal J})))= \left\{\left(\begin{array}{cccc|cccc}
         A&  &       & 0 & B&  &        &0  \\
          &A &       &  &  & B&        &   \\
          &  &\ddots &  &  &  &\ddots  &   \\
         0&  &       & A&  0&  &       &-B \\
         \hline
        -B&  &       &  0& A&  &        & 0 \\
          &-B&       &  &  & A&        &   \\
          &  &\ddots &  &  &  &\ddots  &   \\
          0&  &       & B&  0&  &        &A
         \end{array}\right)\right\}
$$
By Lemma 2.1,  $\theta\circ\psi\circ\theta^{-1}$ (an automorphism
of the algebra of the form (1)) can be extended to an automorphism
of $H(F_n)$. Finally, by interchanging the $k$-th and
$(\frac{n}{2}+k)$-th columns, and $k$-th and $(\frac{n}{2}+k)$-th
rows, $1\le k\le \frac{n}{2}$, and the columns and rows inside the
block (if necessary), we can achieve the following block-diagonal
form:
$$ \left\{\left(\begin{array}{ccccc}
         A& B&        &   & \\
        -B& A&        &  0 & \\
          &  &\ddots  &   & \\
          & 0 &        & A & B\\
          &  &        & -B & A
          \end{array}\right)\right\}\eqno(9)
$$
As a result any subalgebra of $H(F_n)$ of the type $F_m^{(+)}$ can
be brought to the canonical form (2.1). This canonical form is
obviously unique.

{\bf 2.2}. Let $\cal J$ be of the type $H(F_m)$. Acting in the
same manner as before, $\cal J$ can be brought to the unique
canonical form as follows
$$\theta(\psi(\theta^{-1}({\cal J})))=\left\{ \left(\begin{array}{cccccc}
         X&  \ldots   &      0 &         0 &  \ldots &  0 \\
\vdots    &    \ddots &  \vdots&  \vdots   &  \ddots &  \vdots \\
   0      &  \ldots   &      X &         0 &  \ldots &  0  \\
   0      & \ldots    &      0 &         X &  \ldots &  0   \\
\vdots    & \ddots    & \vdots &   \vdots  &  \ddots &  \vdots  \\
   0      & \ldots    &      0 &         0 &  \ldots &       X
    \end{array}\right)\right\},\eqno(10)
$$
where $X$ is a symmetric matrix of order $m$.

{\bf 2.3}. Let $\cal J$ be of the type $H(F_{2m},j)$, $n=4kl$.
Like in the previous cases, by an appropriate automorphism $\psi$
of $F_{\frac{n}{2}}^{(+)}$, $\theta^{-1}(\cal J)$ can be brought
to the following block-diagonal form:

$$\psi(\theta^{-1}({\cal J}))=\left\{ \left(\begin{array}{cccccc}
         X&  \ldots   &      0 &         0 &  \ldots &  0 \\
\vdots    &    \ddots &  \vdots&  \vdots   &  \ddots &  \vdots \\
   0      &  \ldots   &      X &         0 &  \ldots &  0  \\
   0      & \ldots    &      0 &         X &  \ldots &  0   \\
\vdots    & \ddots    & \vdots &   \vdots  &  \ddots &  \vdots  \\
   0      & \ldots    &      0 &         0 &  \ldots &       X
    \end{array}\right)\right\},
$$
where $X$ is a symplectic matrix of order $2m$. If we represent
$X$ as the sum of symmetric and skew-symmetric matrices as
follows:
$$ X=\left(\begin{array}{cc}
               A& -B\\
               B& A
          \end{array}\right)+
\left(\begin{array}{cc}
               -C& -D\\
               -D& C
           \end{array}\right)
$$
where all matrices have order $m$; $A$ is symmetric, $B$, $C$,$D$
are skew-symmetric, then $\theta$ induces the following
representation of $\cal J$ in $H(F_n)$
$$ \theta(\psi(\theta^{-1}({\cal J})))=\left\{\left(\begin{array}{ccc|ccc}
              A& -B&    &   -C& -D& \\
              B&  A&    &   -D&  C& \\
               &   &\ddots   &     &   &\ddots \\
              \hline
              C& D &        &    A& -B& \\
              D&-C &        &    B&  A& \\
               &   & \ddots &     &   &\ddots
               \end{array}\right)\right\}
$$
Similarly, by Lemma 2.1, $\theta\circ\psi\circ\theta^{-1}$ (an
automorphism of the algebra of the form (1)) can be extended to an
automorphism of $H(F_n)$.

By interchanging appropriate blocks, we can reduce it to the
canonical form:
$$   \left\{\left(\begin{array}{ccccccccc}
             A& -B& -C& D&          &  & &  &  \\
             B&  A&  D& C&          &  & &  & \\
             C&  D&  A&-B&          &  &0 &  & \\
             D& -C&  B& A&          &  & &  &\\
              &   &   &  &\ddots    &  & &  & \\
              &   &   &  &         A&-B&-C& D\\
              &   &  0 &  &         B& A& D& C\\
              &   &   &  &         C& D& A&-B\\
              &   &   &  &         D&-C& B& A
              \end{array}\right)\right\}\eqno(11)
$$
From the explicit form (11), the canonical form of $\cal J$ of the
type $H(F_{2m},j)$ is uniquely determined.

\par\medskip
 \par\medskip
If $1\notin {\cal J}$, then $\mbox{rk}\,(e)=k<n$ where $e$ is the
identity element of $\cal J$, and by an appropriate automorphism
of $H(F_n)$ $\cal J$ can be brought to the form:
$$     \left\{\left(\begin{array}{c|c}
                    X& 0 \\
                    \hline
                    0& 0
            \end{array}\right)\right\}\eqno(12)
$$
where $X$ is a symmetric matrix of order $k$. Let $\pi$ denote a
mapping that sends each matrix of the form (12) to the first block
of order $k$. Clearly, ${\cal J}\cong\pi({\cal J})$. Acting in the
same manner as before we can bring $\pi({\cal J})$ to the
canonical form in $H(F_k)$. As a result, the original subalgebra
$\cal J$ also takes the unique canonical form.

\par\medskip

{\bf Case 3}  ${\cal A}=H(F_{2n},j)$

Since the proof of this case is not much different from  the proof
of the case of $H(F_n)$, we will omit most details. According to
classification (see \cite{R2}), any maximal  subalgebra $M$ in
$H(F_{2n},j)$ is isomorphic to one of the following:

1. $H(F_{2k},j)\oplus H(F_{2l},j)$, $k+l=n$,

2. $H(F_{2k},j)\oplus F_l^{(+)}\oplus R$, $k+l=n$. If $k=0$, then
$M=F_n^{(+)}\oplus R$.

3. $J(f,1)$ only if $n=2^{m}$ and either $\dim\,J(f,1)=2(m+1)$,
$m$ is even, or $\dim\,J(f,1)=2m+1$, $m$ is odd.

First we assume that $\cal J$ is a simple matrix subalgebra of
$H(F_{2n},j)$ such that $1\in {\cal J}$. Let $M$ be a maximal
subalgebra which contains $\cal J$, ${\cal J}\subseteq M$. By
Lemma 2.3, ${\cal J}\subset S$. If $S$ is a non-simple semisimple
algebra, then $\cal J$ can be projected into the simple components
of $S$. Hence, the problem will be reduced to the case of
symplectic matrices of order less than $2n$. This reduction stops
only when  either the image  of $\cal J$ can be covered by the
maximal subalgebra with a simple Wedderburn factor $S$ or the
image of $\cal J$  coincides with one of the simple components of
$S$.

 Next we
look into the case when ${\cal J}\subset M$, where $M$ has a
simple Wedderburn factor $S$. There is no loss in generality if we
assume that $M=S\oplus R$, $S\cong F_n^{(+)}$. By Lemma 2.4, $S$
can be brought to the form (5). Notice that any automorphism of
$F_n^{(+)}$ of the form $\varphi(X)=C^{-1}XC$ can be extended to
an automorphism of $H(F_{2n},j)$ as follows:
$$\bar\varphi(X)={\bar C}^{-1}X{\bar C},\quad {\bar C}= \left(\begin{array}{cc}
             C& 0\\
             0& (C^{-1})^t
          \end{array} \right)\eqno(13)
$$

{\bf 3.1} If ${\cal J}\cong H(F_m)$, $m\le n$, then acting by some
automorphism of the form (13), it can be reduced to (3.2). This
canonical form is obviously uniquely determined.

{\bf 3.2} If ${\cal J}\cong F_m^{(+)}$, $m\le n$, then by some
automorphism of the form (13) it can be brought to
$$   \left\{\left(\begin{array}{ccc|ccc}
            X&       &   &   &      & \\
             &\ddots &   &   &  0    &\\
             &       & X^t &   &      &\\
             \hline
             &       &   & X^t &      &   \\
             &   0   &   &   &\ddots& \\
             &       &   &   &      & X
             \end{array}\right)\right\}\eqno (14)
$$
where $X$ is an arbitrary matrix of order $m$. This is the
canonical form (3.1). For further purposes, we introduce the
following characteristic of the canonical form. Let $l$ and $k$
denote the number of $X$-blocks and $X^t$-blocks in the {\it top
left} block of order $n$ in (14). Set $k_{\cal A}({\cal
J})=|l-k|$.With some effort it can be shown that in this case  the
canonical form is also unique. All we have to show is that any two
canonical forms ${\cal J}_1$ and ${\cal J}_2$ of the same type
with $k_{\cal A}({\cal J}_1)\ne k_{\cal A}({\cal J}_2)$ are not
conjugate under symplectic automorphism, or, equivalently,
automorphism of $H(F_{2n},j)$. For clarity, let ${\cal
J}_1=\mbox{diag}\,\{\underbrace{X,X^t,\ldots,X}_{n},X^t,X,\ldots
X^t\}$ and ${\cal J}_2=\mbox{diag}\,\{\underbrace{Y,Y,\ldots
Y^t}_{n},Y^t,Y^t,\ldots, Y\}$ where $X$ and $Y$ are any matrices
of order $m$. Let $\cal S$ stand for the subalgebra of
$H(F_{2n},j)$ of the form (3). Obviously, ${\cal J}_1\subseteq S$,
${\cal J}_2\subseteq S$.
 Next we are going to show that for any automorphism
$\varphi$ of $H(F_{2n},j)$ such that $\varphi({\cal J}_1)={\cal
J}_2$ we can always find a symplectic automorphism $\psi$ that can
be restricted to $\cal S$ and $\psi({\cal J}_1)={\cal J}_2$.
 Let $C$ be a non-singular matrix that  determines $\varphi$. Then, for any $A\in
{\cal J}_1$ there exists $B\in {\cal J}_2$ such that
$$C^{-1}AC=B,\quad AC=CB.\eqno(15)$$
Set $C=(C_{ij})_{i,j=1,s}$ where $C_{ij}$ is a square matrix of
order $m$. By performing a matrix multiplication in (15) we obtain
a series of equations:
$$ XC_{ij}=C_{ij}Y,\quad X^tC_{kl}=C_{kl}Y$$
where $(i,j),(k,l)\in I\times I$, $I=\{1,\ldots,s\}$. Since $X$
and $Y$ can be any matrices of order $m$, $C_{ij}$ can not be
degenerate. Therefore, $Y=C_{ij}^{-1}XC_{ij}$,
$Y=C_{kl}^{-1}X^tC_{kl}$. Hence, the matrix $\bar
C=\mbox{diag}\,\{\underbrace {C_{ij},C_{kl},\ldots,C_{ij}}_n,
 (C_{ij}^t)^{-1}, (C_{kl}^t)^{-1},\ldots, (C_{ij}^t)^{-1}\}$
determines  an automorphism $\psi$ of $H(F_{2n},j)$ such that
$\psi({\cal J}_1)={\cal J}_2$. Besides $\psi$ can be restricted to
$\cal S$, therefore, induces an automorphism of a subalgebra of
the type $F_n^{(+)}$. However we have already showed (case 1.1)
that the two canonical forms in $F_n^{(+)}$ with  $k_{\cal
A}({\cal J}_1)\ne k_{\cal A}({\cal J}_2)$ are not conjugate.

{\bf 3.3} If ${\cal J}\cong H(F_{2m},j)$, $m\le n$, then , it can
be reduced to (3.3).
 Let $l$ and
$k$ denote the number of $X$-blocks and $X^t$-blocks,
respectively, in $A$.  Set $k_{\cal A}({\cal J})=|l-k|$.Then, we
are going to show that any two canonical forms ${\cal J}_1$ and
${\cal J}_2$ of the same type with $k_{\cal A}({\cal J}_1)\ne
k_{\cal A}({\cal J}_2)$ are not conjugate under the automorphism
of $H(F_{2m},j)$. Assume  the contrary, that is, there exists an
automorphism of $H(F_{2m},j)$ such that $\varphi({\cal J}_1)={\cal
J}_2$. Next we can choose ${\cal S}_1\subseteq {\cal J}_1$, ${\cal
S}_1\cong F_m^{(+)}$ such that  $k_{\cal A}({\cal S}_1)=k_{\cal
A}({\cal J}_1)$. Similarly, we can select ${\cal S}_2\subseteq
{\cal J}_2$, ${\cal S}_2\cong F_m^{(+)}$ such that $k_{\cal
A}({\cal S}_2)=k_{\cal A}({\cal J}_2)$. By Lemma 2.2 there exists
$\psi$, $\psi:{\cal J}_2\to {\cal J}_2$, $\psi(\varphi({\cal
S}_1))={\cal S}_2$. From the explicit from of ${\cal J}_2$, $\psi$
can be extended to an automorphism of $H(F_{2m},j)$,
$\psi\circ\varphi:H(F_{2m},j)\to H(F_{2m},j)$. It follows that
${\cal S}_1$ and ${\cal S}_2$ have the same canonical forms, in
particular, $k_{\cal A}({\cal S}_1)=k_{\cal A}({\cal S}_2)$, a
contradiction.

If $1\notin {\cal J}$, then in order to find the canonical form of
$\cal J$ we use the same approach as in the case of $H(F_n)$.

The Theorem is proved.

\end{proof}
\par\medskip
\par\medskip

Let $e$ denote the identity element of ${\cal J}\cong F_m^{(+)}$,
and $\rho$ stand for the natural representation of $\cal J$ in
$F^n$, $m<n$. Obviously, $\rho$ induces the representation of
$S({\cal J})$ in $F^n$. If $S({\cal J})$ is a non-simple
semisimple associative algebra, that is, $S({\cal J})={\cal
I}_1\oplus {\cal I}_2$ where ${\cal I}_1$, ${\cal I}_2$ are
isomorphic simple ideals, then $\rho=\rho_1\oplus \rho_2$ where
$\rho_i$ is a representation of ${\cal I}_i$ in the corresponding
invariant subspace of $F^n$. Then $k_{\cal A}({\cal
J})=|\mbox{deg}\,\rho_1(I_1)-\mbox{deg}\,\rho_2(I_2)|$.

\begin{theorem}
Let $\cal A$  be a Jordan algebra of any of the following types:
$F_n^{(+)}$, $H(F_n)$ or $H(F_{2n},j)$, $n\ge 3$, and $\cal J$,
${\cal J}'$ be proper simple matrix subalgebras of $\cal A$. If
${\cal J}'$  has the same type as $\cal J$ does, then ${\cal
J}'\in C({\cal J})$ if and only if

 (1) $\mbox{rk}\,(e)= \mbox{rk}\,(e')$;

 (2) $k_{\cal A}({\cal J})= k_{\cal A}({\cal J}'),$
 in the case when ${\cal J}\cong F_m^{(+)}$ or $H(F_{2m},j)$, for some $m<n$, and
 ${\cal A}\cong F_n^{(+)}$ or $H(F_{2n},j)$.

\end{theorem}

\begin{proof}
First it should be noted that the degree of ${\cal J}\ge 3$. The
case of $\cal J$ of the degree 2 will be considered later in the
text.

{\bf The case of $F_n^{(+)}$}

In this case we assume that $\cal J$ and ${\cal J}'$ are
subalgebras of $F_n^{(+)}$ which is as usual the set of all
matrices of order $n$ closed under the Jordan multiplication. This
case breaks into the following subcases.

{\bf (1)} Let $\cal J$ be of the type $F_m^{(+)}$ for some $m<n$.
 First we assume that $S({\cal J})$ is a simple algebra.
Equivalently, $k_{\cal A}({\cal J})=\mbox{rk}(e)$. Let ${\cal J}'$
be as given in the conditions of the Theorem. If ${\cal J}'\in
{\cal C}({\cal J})$, then there exists an automorphism $\varphi$
of $F_n^{(+)}$ which maps ${\cal J}'$ onto $\cal J$. It follows
that $\varphi(e')=e,$ therefore, $\mbox{rk}\,(e')=\mbox{rk}\,(e)$.
Besides, $\varphi(S({\cal J}'))=S({\cal J})$. Hence, $S({\cal
J}')$ is also simple, $k({\cal J}')=\mbox{rk}(e')$.  It follows
that $k({\cal J}')=\mbox{rk}\,(e')=\mbox{rk}\,(e)=k({\cal J})$.

Conversely, if $\mbox{rk}(e')= \mbox{rk}(e)$ and  $k({\cal J})=
k({\cal J}')$, then $k({\cal J}')= k({\cal J})=\mbox{rk}(e)=
\mbox{rk}(e')$, because $k({\cal J})=\mbox{rk}(e)$. Therefore,
$k({\cal J}')=\mbox{rk}(e')$, that is, $S({\cal J}')$ is also
simple, and $\cal J$, ${\cal J}'$ have the same canonical forms.
This implies that ${\cal J}'\in {\cal C}({\cal J})$.

Now we assume that $S({\cal J})$ is a non-simple semisimple
subalgebra. Let ${\cal J}'$ be another subalgebra which satisfies
the conditions of the Theorem. If ${\cal J}'\in {\cal C}({\cal
J})$, then there exists an automorphism $\varphi$ of $F_n^{(+)}$
which maps ${\cal J}'$ onto $\cal J$. Therefore, ${\cal J}'$ and
$\cal J$ have equivalent representations in $F^n$, and so do
$S({\cal J}')$ and $S({\cal J})$. Consequently, either
$\mbox{deg}\,\rho_1({\cal I}_1)=\mbox{deg}\,\rho_1({\cal I}_1')$
and $\mbox{deg}\,\rho_2({\cal I}_2)=\mbox{deg}\,\rho_2({\cal
I}_2')$ or $\mbox{deg}\,\rho_1({\cal
I}_1)=\mbox{deg}\,\rho_2({\cal I}_2')$ and
$\mbox{deg}\,\rho_2({\cal I}_2)=\mbox{deg}\,\rho_1({\cal I}_1')$.
Equivalently, $|\mbox{deg}\rho_1({\cal
I}_1)-\mbox{deg}\rho_2({\cal I}_2)|=|\mbox{deg}\rho_1({\cal
I}'_1)-\mbox{deg}\rho_2({\cal I}'_2)|$, that is, $k_{\cal A}({\cal
J})=k_{\cal A}({\cal J}')$.

Conversely, if  $\mbox{rk}(e')= \mbox{rk}(e)$ and $k_{\cal
A}({\cal J})=k_{\cal A}({\cal J}')$, then $\cal J$ and ${\cal J}'$
have the same canonical forms. Therefore, these subalgebras are
conjugate by some automorphism of $F_n^{(+)}$, and ${\cal J}'\in
{\cal C}({\cal J})$.

{\bf (2)} Let $\cal J$ be of the type $H(F_m)$ for some $m\le n$.

Suppose that ${\cal J}'$ is another subalgebra of $F_n^{(+)}$
which has the type $H(F_m)$. If ${\cal J}'$ is conjugate to $\cal
J$ under some automorphism $\varphi$ of $F_n^{(+)}$ then
$\varphi(e')=e$ and $\mbox{rk}(e')=\mbox{rk}(e)$. In other words,
the canonical form of ${\cal J}'$ is exactly the same as that of
$\cal J$. Conversely, if $\mbox{rk}(e')=\mbox{rk}(e)$, then $\cal
J$ and ${\cal J}'$ have the same canonical forms. Therefore,
${\cal J}'\in {\cal C}({\cal J})$.

{\bf (3)} Let $\cal J$ be of the type $H(F_{2m},j)$ for some $m\le
n$. The proof of this case  is exactly the same as the previous
proof.

{\bf The  case of $H(F_n)$}

Suppose that $\cal J$ and ${\cal J}'$ are two subalgebras of
$H(F_n)$ that satisfy the conditions of the Theorem.

{\bf (1)} Let $\cal J$ as well as ${\cal J}'$ be of the type
$F_m^{(+)}$. Assume that ${\cal J}'\in {\cal C}({\cal J})$. It
follows that there exists an automorphism of $H(F_n)$ such that
$\varphi({\cal J}')={\cal J}$. Hence,
$\mbox{rk}(e')=\mbox{rk}(e)$. Conversely, if
$\mbox{rk}(e')=\mbox{rk}(e)$, then $\cal J$ and ${\cal J}'$ have
the same canonical form. Therefore, ${\cal J}'\in {\cal C}({\cal
J})$.

{\bf (2)} Now let both $\cal J$ and ${\cal J}'$ have the type
$H(F_k)$ (or $H(F_{2k},j)$). If ${\cal J}'\in {\cal C}({\cal J})$,
then there exists an automorphism $\varphi$ of $H(F_n)$  that
sends ${\cal J}'$ onto $\cal J$, $\varphi({\cal J}')={\cal J}$.
Consequently, $\mbox{rk}(e')=\mbox{rk}(e)$.

Conversely, if $\mbox{rk}(e')=\mbox{rk}(e)$, then they have the
same canonical form. Therefore, ${\cal J}'\in {\cal C}({\cal J})$.

\par\medskip
{\bf The case of $H(F_{2n},j)$}

 Suppose that $\cal J$ and ${\cal
J}'$ are two subalgebras of $H(F_{2n},j)$ that satisfy the
conditions of the Theorem.

{\bf (1)} Let $\cal J$ as well as ${\cal J}'$ be of the type
$F_m^{(+)}$, $m<n$. Assume that ${\cal J}'\in {\cal C}({\cal J})$.
It follows that there exists an automorphism of $H(F_{2n},j)$ such
that $\varphi({\cal J}')={\cal J}$. Hence,
$\mbox{rk}(e')=\mbox{rk}(e)$. Since ${\cal J}'$ and $\cal J$ are
conjugate in $H(F_{2n},j)$, they  have the same canonical forms in
$H(F_{2n},j)$. Therefore, $k_{\cal A}({\cal J})= k_{\cal A}({\cal
J}')$.

Conversely, if all conditions hold true, then $\cal J$ and ${\cal
J}'$ have the same canonical forms. Therefore, ${\cal J}'\in {\cal
C}({\cal J})$.

{\bf (2)} Now let both $\cal J$ and ${\cal J}'$ have the type
$H(F_m)$, $m<n$. If ${\cal J}'\in {\cal C}({\cal J})$, then there
exists an automorphism of $H(F_{2n},j)$ that sends ${\cal J}'$
onto $\cal J$, $\varphi({\cal J}')=\cal J$. Consequently,
$\mbox{rk}(e')=\mbox{rk}(e)$.

Conversely, if $\mbox{rk}(e')=\mbox{rk}(e)$, then they have the
same canonical forms. Therefore, ${\cal J}'\in {\cal C}({\cal
J})$.

{\bf (3)} Now let both $\cal J$ and ${\cal J}'$ have the type
$H(F_{2m},j)$, $m<n$. If ${\cal J}'\in {\cal C}({\cal J})$, then
there exists an automorphism of $H(F_{2n},j)$ that sends ${\cal
J}'$ onto $\cal J$, $\varphi({\cal J}')=\cal J$. Consequently,
$\mbox{rk}(e')=\mbox{rk}(e)$, $k_{\cal A}({\cal J})= k_{\cal
A}({\cal J}')$.

Conversely, if $\mbox{rk}(e')=\mbox{rk}(e)$ and $k_{\cal A}({\cal
J})= k_{\cal A}({\cal J}')$, then they have the same canonical
forms. Therefore, ${\cal J}'\in {\cal C}({\cal J})$

 The Theorem is proved.

\end{proof}

\begin{corollary}
If $m$ is any number such that $m\le n$, and $n=mk+r$, $0\le r<m$,
then there exist subalgebras of $F_n^{(+)}$ of the type $H(F_m)$.
Moreover, there are precisely $k$ conjugate classes corresponding
to $H(F_m)$. If $2m\le n$, and $n=2mk+r$, $0\le r<m$ then
$F_n^{(+)}$ has subalgebras of the type $H(F_{2m},j)$, and the
number of conjugate classes corresponding to $H(F_{2m},j)$ is
equal to $k$. Finally, if $m<n$, and $n=mk+r$, $0\le r<m$ then
there exist subalgebras of $F_n^{(+)}$ of the type $F_m^{(+)}$,
and, moreover, the number of conjugate classes is given by
$\sum_{j=1}^k [\frac{j}{2}].$
\end{corollary}

\begin{corollary}
If $m$ is any number such that $m<n$, and $n=mk+r$, $0\le r<m$,
then there exist subalgebras of $H(F_n)$ of the type $H(F_m)$.
Moreover, there are precisely $k$ conjugate classes corresponding
to $H(F_m)$. If $2m\le n$, and $n=2mk+r$, $0\le r<m$ then $H(F_n)$
has subalgebras of the type $F_m^{(+)}$, and the number of
conjugate classes corresponding to $F_m^{(+)}$ is equal to $k$.
Finally, if $ 4m\le n$, and $n=4mk+r$, $0\le r<m$ then there exist
subalgebras of $H(F_n)$ of the type $H(F_{2m},j)$, and, moreover,
the number of conjugate classes is $k$.
\end{corollary}

\begin{corollary}
If $m$ is any number such that $m\le n$, and $n=mk+r$, $0\le r<m$,
then there exist subalgebras of $H(F_{2n},j)$ of the type
$H(F_m)$. Moreover, there are precisely $k$ conjugate classes
corresponding to $H(F_m)$. If $m\le n$, and $n=mk+r$, $0\le r<m$
then $H(F_{2n},j)$ has subalgebras of the type $F_m^{(+)}$, and
the number of conjugate classes corresponding to $F_m^{(+)}$ is
equal to $\sum_{j=1}^k [\frac{j}{2}].$ Finally, if $ m<n$, and
$n=mk+r$, $0\le r<m$ then there exist subalgebras of $H(F_n,j)$ of
the type $H(F_m,j)$, and, moreover, the number of conjugate
classes is $\sum_{j=1}^k [\frac{j}{2}].$
\end{corollary}

\subsection{Subalgebras of the type $J(f,1)$}

First we recall a few facts concerning Clifford algebras over a
field of characteristic not 2(see \cite{R2}). Let ${\cal
J}=F\oplus V$ where $V=\mbox{span}\langle
x_1,\ldots,x_{2m}\rangle$, and $f$ a non-degenerate symmetric
bilinear form on $V$. Then, $C(V,f)$ is a central simple
associative algebra with a unique canonical involution "---" such
that it fixes elements from $V$. In this case the imbedding of
$\cal J$  into $C(V,f)^{(+)}$ we will call {\it canonical of the
first type}. Next, let ${\cal J}=F\oplus V$ where
$V=\mbox{span}\langle x_1,\ldots,x_{2m+1}\rangle$, and
$V_0=\mbox{span}\langle x_1,\ldots,x_{2m}\rangle$.  Then, $C(V,f)$
 is isomorphic to a tensor product of
$C(V_0,f)$  and the two-dimensional center $E$ of $C(V,f)$.
Moreover, $E=F[z]$ where $z=x_1x_2\ldots x_{2m+1}$. In other
words, $C(V,f)={\cal I}_1\oplus {\cal I}_2$, ${\cal I}_i\cong
C(V_0,f)$. Note that $F\oplus V\cong {\cal J}/{\cal I}_i\subseteq
C(V,f)/{\cal I}_i\cong C(V_0,f)$. This imbedding of ${\cal
J}=F\oplus V$ into $C(V_0,f)^{(+)}$ we will call {\it canonical of
the second type.}

 Let $\cal A$ be a simple matrix Jordan algebra, and $\cal J$ be a
 subalgebra of the type $J(f,1)$. According to \cite{R2}, $\cal J$
 of the type $J(f,1)$ is  maximal in $\cal A$ if and only if one
 of the following cases hold
\par\medskip

{\bf 1}. ${\cal A}=(C(V_0,f),-)$, ${\cal J}=F\oplus V$ where
$\dim V=2m+1$ and  $m$ is odd.

{\bf 2}. ${\cal A}=H(C(V_0,f),-)$, ${\cal J}=F\oplus V$ where
$\dim V=2m+1$,  $m$ is even.

{\bf 3}. ${\cal A}=H(C(V,f),-)$, ${\cal J}=F\oplus V$ where $\dim
V=2m$
\par\medskip
Next we recall that if $\dim\,V=2m$, and $m\equiv 0,1(\mbox{mod}\,
4)$ then $\dim\, H(C(V,f),-)=2^{m-1}(2^m+1)$. If $\dim\,V=2m$ and
$m\equiv 2,3(\mbox{mod}\, 4)$ then  $\dim\,H(C(V,f),-)=
2^{m-1}(2^m-1)$. If $\dim\,V=2m+1$ and $m\equiv 0(\mbox{mod}\, 4)$
then $\dim\,H(C(V_0,f),-)=2^{m-1}(2^m+1)$. If $\dim\,V=2m+1$ and
$m\equiv 2(\mbox{mod}\, 4)$ then
$\dim\,H(C(V_0,f),-)=2^{m-1}(2^m-1)$.
\par\medskip
\par\medskip
{\bf\large  Canonical realizations of $J(f,1)$}
\par\medskip
Let $\cal A$ be a simple matrix Jordan algebra, and ${\cal
J}=F\oplus V$ is a subalgebra of $\cal A$.

\noindent 1.1 ${\cal A}=F_n^{(+)}$, $n=2^ml+r$, $\dim\,V=2m$,
$${\cal J}=\{\mbox{diag}(\underbrace{X,\ldots,X}_l ,0,\ldots,0)\}$$
where $X$ is a matrix of order $2^m$, and  if $\pi_i$ denotes the
projection on the $i$th non-zero block, then   $\pi_i({\cal
J})\subseteq F_{2^m}^{(+)}$ is a canonical imbedding of the first
type.
\par\medskip
\noindent 1.2 ${\cal A}=F_n^{(+)}$, $n=2^ml+r$, $\dim\,V=2m+1$,
$${\cal J}=\{\mbox{diag}(\underbrace{X,\ldots,X}_l ,0,\ldots,0)\}$$
where  $X$ is a matrix of order $2^m$, and  $\pi_i({\cal
J})\subseteq F_{2^m}^{(+)}$ is a canonical imbedding of the second
type.
\par\medskip
\noindent 1.3 ${\cal A}=F_n^{(+)}$, $n=2^ml+r$, $\dim\,V=2m+1$,
$${\cal
J}=\{\mbox{diag}(\underbrace{X,\ldots,X,}_s\underbrace{X^t,\ldots,X^t}_k
0,\ldots,0)\}$$ where $s+k=l$, $X$ is a matrix of order $2^m$, and
$\pi_i({\cal J})\subseteq F_{2^m}^{(+)}$ is a canonical imbedding
of the second type.
\par\medskip
\noindent 2.1 ${\cal A}=H(F_n)$, $n=2^{m}l+r$, $\dim\,V=2m$,
$${\cal J}=\{\mbox{diag}(\underbrace{X,\ldots,X}_l ,0,\ldots,0)\}$$
where  $X$ is a symmetric matrix of order $2^m$, and $\pi_i({\cal
J})\subseteq F_{2^m}^{(+)}$ is a canonical imbedding of the first
type.
\par\medskip
\noindent 2.2 ${\cal A}=H(F_n)$, $n=2^{m+1}l+r$, $\dim\,V=2m$,
$${\cal J}=\{\mbox{diag}(\underbrace{X,\ldots,X}_l ,0,\ldots,0)\}$$
where $X$ is of the form (1) in which $A$ and $B$ are of order
$2^m$. If $\cal S$ denotes the algebra of the form (1), then
$\pi_i({\cal J})\subseteq {\cal S}$  is a canonical imbedding of
the first type.
\par\medskip
\noindent 2.3 ${\cal A}=H(F_n)$, $n=2^{m+1}l+r$, $\dim\,V=2m+1$,
$${\cal J}=\{\mbox{diag}(\underbrace{X,\ldots,X}_l ,0,\ldots,0)\}$$
where $X$ is of the form (1) in which $A$ and $B$ are of order
$2^m$. If $\cal S$ denotes the  entire algebra of the form (1),
then $\pi_i({\cal J})\subseteq {\cal S}$  is a canonical imbedding
of the second type.
\par\medskip
\noindent 2.4 ${\cal A}=H(F_n)$, $n=2^ml+r$, $\dim\,V=2m+1$,
$${\cal J}=\{\mbox{diag}(\underbrace{X,\ldots,X}_l ,0,\ldots,0)\}$$
where $X$ is a symmetric matrix, and   $\pi_i({\cal J})\subseteq
F_{2^m}^{(+)}$ is a canonical imbedding of the second type.
\par\medskip
\noindent 3.1 ${\cal A}=H(F_{2n},j)$, $n=2^ml+r$, $\dim\,V=2m$,
$${\cal J}=\{\mbox{diag}(\underbrace{X,\ldots,X}_l ,\underbrace{0,\ldots,0}_k,
\underbrace{X,\ldots,X}_l ,\underbrace{0,\ldots,0)}_k,\}$$ where
$k+l=n$, $X$ is a symmetric matrix of order $2^m$, and
$\pi_i({\cal J})\subseteq F_{2^m}^{(+)}$ is a canonical imbedding
of the first type.

\noindent 3.2 ${\cal A}=H(F_{2n},j)$, $n=2^ml+r$, $\dim\,V=2m$,
$\cal J$ has a canonical form (3.3), and if $\pi_i$ denotes the
projection of $\cal J$ into $i$th simple component (of the type
$H(F_{2^m},j)$) of (3.3), then $\pi_i({\cal J})\subseteq
H(F_{2^m},j)$ is a canonical imbedding of the first type.
\par\medskip
\noindent 3.3 ${\cal A}=H(F_{2n},j)$, $n=2^ml+r$, $\dim\,V=2m+1$,
$\cal J$ has a canonical form (3.3) where $\pi_i({\cal
J})\subseteq H(F_{2^m},j)$ is a canonical imbedding of the second
type.
\par\medskip
\noindent 3.4 ${\cal A}=H(F_{2n},j)$, $n=2^ml+r$, $\dim\,V=2m+1$,
$${\cal
J}=\{\mbox{diag}(\underbrace{X,\ldots,X,}_s\underbrace{X^t,\ldots,X^t}_k
0,\ldots,0,\underbrace{X^t,\ldots,X^t,}_s\underbrace{X,\ldots,X}_k
0,\ldots,0 )\}$$ where $s+k=l$, $X$ is a matrix of order $2^m$,
and $\pi_i({\cal J})\subseteq F_{2^m}^{(+)}$ is a canonical
imbedding of the second type.

\par\medskip
\begin{theorem}
Let $\cal A$ be a simple matrix Jordan algebra, and $\cal J$ be a
subalgebra of $\cal A$ of the type $J(f,1)$. Then, $\cal J$ has a
unique canonical form as above.
\end{theorem}
\begin{proof}
Let ${\cal J}=F\oplus V$. Then the following cases occur.

{\bf Case} ${\cal A}=F_n^{(+)}$

{\bf 1.1} Let $\dim\, V=2m$. Then $U({\cal J})\cong C(V,f)$ is a
simple algebra. In particular, $S({\cal J})\cong U({\cal J})$.

If $S({\cal J})={\cal A}$, then $n=2^m$, ${\cal A}\cong U({\cal
J})$.  Therefore, the imbedding of $\cal J$ into $\cal A$ is
equivalent to the imbedding of $F\oplus V$ into $C(V,f)^{(+)}$.
Therefore, this is a canonical imbedding of the first type.

If $S({\cal J})\subset {\cal A}$, then $S({\cal J})$ is a proper
simple associative subalgebra of $F_n$. Therefore, $S({\cal J})$
can be reduced to
$$ \{\mbox{diag}(\underbrace{Y,\ldots,Y}_l, 0,\ldots,0)\}\eqno (16)$$
where the order of $Y$ is $2^m$, and $n=2^ml+r$. As a result,
$\cal J$ also takes the canonical form 1.1.

{\bf 1.2} Let $\dim\,V=2m+1$. Then  $U({\cal J})\cong C(V,f)$, and
$U({\cal J})={\cal I}_1\oplus{\cal I}_2$, ${\cal I}_i\cong
C(V_0,f)$. Hence $S({\cal J})$ is isomorphic to either $C(V,f)$ or
$C(V_0,f)$.

If $S({\cal J})={\cal A}$, then the imbedding of $\cal J$ into
$\cal A$ is the canonical of the second type.

If $S({\cal J})\cong {\cal I}_i$, and $S({\cal J})\subset {\cal
A}$, then $S({\cal J})$ is a proper simple associative subalgebra
of $F_n$. Therefore, $S({\cal J})$ can be reduced to (16).  As a
result, $\cal J$ takes the canonical form 1.2.

Finally, if $S({\cal J})={\cal I}_1\oplus{\cal I}_2$, then $\cal
J$ takes the canonical form 1.3

{\bf Case} ${\cal A}=H(F_n)$

Let $M$ be the maximal subalgebra of $H(F_n)$ such that ${\cal
J}\subseteq M\subseteq H(F_n)$. Then, the following cases occur.

1. $M=S\oplus R$ where $S=S_1\oplus S_2$ a semisimple factor, $R$
the radical. Then, we reduce the problem to the case of symmetric
matrices of a lower dimension (see section 2.1).

2. $M=S$ where $S=S_1\oplus S_2$. Like in the previous case we can
reduce the problem to the case of symmetric matrices of a lower
dimension.

3. $M=S\oplus R$ where $S\cong F_{\frac{n}{2}}^{(+)}$, $R$ the
radical.

4. $M=F\oplus W$ where $W$ is a finite-dimensional vector space.

After a series of reductions of the form 1 and 2, the image of
$\cal J$ becomes a subalgebra of

$$\left(\begin{array}{ccccc}
{\cal A}_1&             &             &    0    &   \\
          & \ddots      &             &        & \\
          &             & {\cal A}_i  &         & \\
          &             &             &\ddots   &\\
         &  0           &             &        & {\cal A}_k
          \end{array}\right)
$$
where ${\cal A}_i\cong H(F_{n_i})$. Let $\pi_i$ be the projection
of $\cal J$ into  ${\cal A}_i$. To simplify our notations we
denote $\pi_i({\cal J})$ as ${\cal J}'$, and the maximal
subalgebra of ${\cal A}_i$ which covers ${\cal J}'$ as $M_i$,
${\cal J}'\subseteq M_i\subseteq {\cal A}_i=H(F_{n_i})$.

{\bf Case 1.} Let $\dim\, V=2m$ and $m\equiv 0,1 (\mbox{mod} 4)$.
Then we have the following cases:

(a) Let $M_i=F\oplus W$. If $S({\cal J}')=F_{n_i}$, then
$n_i=2^m$, $F_{n_i}\cong C(V,f)$, and the imbedding of ${\cal J}'$
into $F_{n_i}^{(+)}$ is equivalent to the imbedding of $F\oplus V$
into $C(V,f)^{(+)}$, that is, canonical of the first type.  If
$S({\cal J}')\subset F_{n_i}$, then $H(S({\cal J}'))\subseteq
H(F_{n_i})$ is a proper subalgebra of $H(F_{n_i})$. Hence,
$n_i=2^ml+r$, and $H(S({\cal J}'))$ can be reduced to (16) in
which $X$ denotes a symmetric matrix of order $2^m$. Then, ${\cal
J}'$ takes the canonical form 2.1.

(b) Let $M_i=S\oplus R$ where $S\cong F_{\frac{n_i}{2}}^{(+)}$. By
using $\theta$-isomorphism (see section 2.1) we obtain that
$\theta^{-1}({\cal J}')\subseteq F_{\frac{n_i}{2}}^{(+)}$. If
$S(\theta^{-1}({\cal J}'))=F_{\frac{n_i}{2}}^{(+)}$, then
$n_i=2^{m+1}$ and the imbedding of $\theta^{-1}({\cal J}')$ into
$F_{\frac{n_i}{2}}^{(+)}$ is the canonical imbedding of the first
type. In particular, $\theta^{-1}({\cal J}')\subseteq
H(F_{\frac{n}{2}})$. As a result ${\cal J}'$ takes the canonical
form 2.1 in which $l=2$ and no zeros. If $S({\cal J}')\subset
F_{\frac{n_i}{2}}^{(+)}$, then $S(\theta^{-1}({\cal J}'))$ is a
proper simple subalgebra of $ F_{\frac{n_i}{2}}$, therefore, takes
the  form (16) and $n_i=2^{m+1}l+r$. Hence ${\cal J}$ takes the
canonical form 2.1.

{\bf Case 2} Let $\dim\, V=2m$, $m\equiv 2,3 (\mbox{mod} 4)$.

(a) Let $M_i=F\oplus W$. If $S({\cal J}')=F_{n_i}$, then
$n_i=2^m$, $F_{n_i}\cong C(V,f)$, ${\cal J}'\subseteq H(F_{n_i})$.
Hence we have the following commutative diagram:

$$             \begin{array}{ccc}
             {\cal J}'=F\oplus V& \stackrel{id}{\longrightarrow}& {\cal J}'=F\oplus V\\
             \downarrow\lefteqn{\sigma}&    & \downarrow\lefteqn{\eta}\\
             U({\cal J}')&  \stackrel{\varphi}{\longleftarrow} &
             F_{n_i}^{(+)}
              \end{array}
$$
where $\sigma=\varphi\circ\eta$. Therefore, $\sigma({\cal
J}')=\varphi(\eta({\cal J}'))$ is symmetric with respect to the
canonical involution "---" which is symplectic in this particular
case. On the other hand, $\sigma({\cal J}')$ is also symmetric
with respect to $j'=\varphi\circ t\circ\varphi^{-1}$. By the
uniqueness of "---", $j'$ equals to "---". However it is not
possible because $\dim\, H(C(V,f),-)=\frac{2^m(2^m-1)}{2}\ne
\frac{2^m(2^m+1)}{2}=\dim\, H(C(V,f),j')$. If $S({\cal J}')\subset
F_{n_i}$, then $H(S({\cal J}'))\subseteq H(F_{n_i})$ is a proper
subalgebra of $H(F_{n_i})$. Hence $n_i=2^ml+r$, and $H(S({\cal
J}'))$ can be reduced to (16) where $X$ denotes a symmetric matrix
of order $2^m$. Let $\pi_{ij}$ denote the projection on $j$th
non-zero block of (16). Then the imbedding $\pi_{ij}({\cal
J}')\subseteq \pi_{ij}(H(S({\cal J}')))$ is similar to the above
imbedding, which is not possible.

(b) Let $M_i=S\oplus R$ where $S\cong F_{\frac{n_i}{2}}^{(+)}$.
Then $\theta^{-1}({\cal J}')\subseteq F_{\frac{n_i}{2}}^{(+)}$.
Since $S(\theta^{-1}({\cal J}'))\cong  U({\cal J}')$, then
$n_i=2^{m+1}l+r$, and $S(\theta^{-1}({\cal J}'))$ can be reduced
to (16) in which $X$ is any matrix of order $2^m$. Hence $\pi_{ij}
(\theta^{-1}({\cal J}'))\subset F_{2^m}^{(+)}$  is a canonical
imbedding of the first type, and ${\cal J}$ has the canonical form
2.2.

{\bf Case 3} Let $\dim\, V=2m+1$ where $m$ is odd.

(a) Let $M_i=F\oplus W$. If $S({\cal J}')=F_{n_i}$, then
$n_i=2^m$, $F_{n_i}\cong C(V_0,f)$. Therefore, the imbedding of
${\cal J}'$ into $F_{n_i}^{(+)}$ is equivalent to the imbedding of
$F\oplus V$ into $C(V_0,f)^{(+)}$ which is canonical imbedding of
the second type. Since $m$ is odd, ${\cal J}'$ is a maximal
subalgebra in $F_{n_i}^{(+)}$. However, ${\cal J}'\subseteq
H(F_{n_i})$, hence, ${\cal J}'$ cannot be maximal. This case is
not possible. If $S({\cal J}')\subseteq F_{n_i}$, then $H(S({\cal
J}'))\subseteq H(F_{n_i})$ is a proper subalgebra of $H(F_{n_i})$,
therefore, can be reduced to (16). However, the imbedding of
$\pi_{ij}({\cal J}')$ into $F_{2^m}^{(+)}$ is as shown above.
Hence this case is also not possible.

(b)Let $M_i=S\oplus R$ where $S\cong F_{\frac{n_i}{2}}^{(+)}$.
Acting in the same manner as in case 2(b) we will come to the
canonical form 2.3.

{\bf Case 4.} Let $\dim\, V=2m+1$ and $m\equiv 0 (\mbox{mod} 4)$.
Acting in the same manner as in previous cases we will reduce
${\cal J}'$ to the canonical form 2.4.

{\bf Case 5.} Let $\dim\, V=2m+1$ and $m\equiv 2 (\mbox{mod} 4)$.
Acting in the same manner as in previous cases we will reduce
${\cal J}'$ to the canonical form 2.3.

{\bf Case} ${\cal A}=H(F_{2n},j)$

Let $M$ be the maximal subalgebra of $H(F_{2n},j)$ such that
${\cal J}\subseteq M\subseteq H(F_{2n},j)$. Then, the following
cases occur.

1. $M=S\oplus R$ where $S=S_1\oplus S_2$ a semisimple factor, $R$
the radical. Then, we reduce the problem to the case of symplectic
matrices of a lower dimension (see section 2.1).

2. $M=S$ where $S=S_1\oplus S_2$. Like in the previous case we can
reduce the problem to the case of symplectic matrices of a lower
dimension.

3. $M=S\oplus R$ where $S\cong F_{n}^{(+)}$, $R$ the radical.

4. $M=F\oplus W$ where $W$ is a finite-dimensional vector space.

After a series of reductions of the form  1 and 2, the image of
$\cal J$ becomes a subalgebra of the algebra in the canonical form
(3.3) in which the $i$th component has order $2n_i$. Let $\pi_i$
denote the projection of $\cal J$ into the $i$th simple component
of (3.3).

{\bf Case 1.} Let $\dim\, V=2m$ and $m\equiv 0,1 (\mbox{mod} 4)$.

(a) Let $M_i=F\oplus W$. If $S({\cal J}')=F_{{2n}_i}$,
$F_{{2n}_i}\cong C(V,f)$, $2n_i=2^m$.  Acting in the same manner
as in case 2(a), we can show that this situation is not possible.
Likewise if $S({\cal J}')\subset F_{{2n}_i}$ then we can reduce
this case to the case just considered. Therefore, it also never
occurs.

(b) Let $M_i=S\oplus R$ where $S\cong F_{n_i}^{(+)}$. Then ${\cal
J}'\subset F_{{n}_i}^{(+)}$, therefore, $S({\cal J}')$ can be
brought to (16), and $\pi_{ij}({\cal J}')\subset F_{2^m}^{(+)}$ is
the canonical imbedding of the first type. Finally the original
subalgebra takes the form 3.1

{\bf Case 2} Let $\dim\, V=2m$, $m\equiv 2,3 (\mbox{mod} 4)$.

(a) Let $M_i=F\oplus W$. If $S({\cal J}')=F_{2n_i}$, then
$2n_i=2^m$, $F_{2n_i}\cong C(V,f)$, ${\cal J}'\subseteq
F_{2n_i}^{(+)}$ is the canonical imbedding of the first type.
 If $S({\cal J}')\subset F_{2n_i}$, then $H(S({\cal J}'),j)\subseteq H(F_{2n_i},j)$
is a proper subalgebra of $H(F_{2n_i},j)$, that is, $n_i=2^ml+r$,
and $H(S({\cal J}'),j)$ can be reduced to (3.3) in which each
component has order $2^m$. Then, ${\cal J}$ takes the canonical
form 3.2.

(b) Let $M_i=S\oplus R$ where $S\cong F_{n_i}^{(+)}$. This case
also lead us to the canonical form 3.2.

{\bf Case 3} Let $\dim\, V=2m+1$ where $m$ is odd.

(a) Let $M_i=F\oplus W$. If $S({\cal J}')=F_{2n_i}$, then
$2n_i=2^m$, $F_{2n_i}\cong C(V_0,f)$. Therefore, the imbedding of
${\cal J}'$ into $F_{2n_i}^{(+)}$ is equivalent to the imbedding
of $F\oplus V$ into $C(V_0,f)^{(+)}$ which is canonical imbedding
of the second type. Since $m$ is odd, ${\cal J}'$ is a maximal
subalgebra in $F_{2n_i}^{(+)}$. However, ${\cal J}'\subseteq
H(F_{2n_i},j)$, hence, ${\cal J}'$ cannot be maximal. This case is
not possible. If $S({\cal J}')\subset F_{2n_i}$, then $H(S({\cal
J}'),j)\subset H(F_{2n_i},j)$ is a proper subalgebra of
$H(F_{2n_i},j)$, therefore, can be reduced to (3.3). Let
$\pi_{ij}$ denote the projection of ${\cal J}'$ into the $j$th
simple component of (3.3). However, the imbedding of
$\pi_{ij}({\cal J}')$ into $F_{2^m}^{(+)}$ is as shown above.
Hence this case is also not possible.

(b) Let $M_i=S\oplus R$ where $S\cong F_{n_i}^{(+)}$. Then ${\cal
J}'\subset F_{{n}_i}^{(+)}$, therefore, $S({\cal J}')$ can be
brought to (16), and $\pi_{ij}({\cal J}')\subset F_{2^m}^{(+)}$ is
the canonical imbedding of the second type. Finally the original
subalgebra takes the form 3.4

{\bf Case 4.} Let $\dim\, V=2m+1$ and $m\equiv 0 (\mbox{mod} 4)$.
Acting in the same manner as in previous cases we will reduce
${\cal J}'$ to the canonical form 3.1.

{\bf Case 5.} Let $\dim\, V=2m+1$ and $m\equiv 2 (\mbox{mod} 4)$.
Acting in the same manner as in previous cases we will reduce
${\cal J}'$ to the canonical form 3.2.

\par\medskip
\begin{corollary}
Let $\cal A$ be a simple matrix Jordan algebra of degree $\ge 3$,
and ${\cal J}=F\oplus V$ where either $\dim\,V=2m$ or
$\dim\,V=2m+1$. Then, $\cal J$ is a subalgebra of $\cal A$ if and
only if

1. $2^m\le n$,

2.  $2^{m+1}\le n$, in the case when ${\cal A}=H(F_n)$ and
$m\equiv 2,3(\mbox{mod}\, 4)$.

\end{corollary}

\end{proof}

\par\medskip
\par\medskip

The author uses this opportunity to thank her supervisor Prof.
Bahturin for his helpful cooperation, many useful  ideas and
suggestions.

\end{document}